\documentclass[12pt,reqno]{amsart} 

\usepackage[arrow,matrix,curve]{xy}

\usepackage[dvips]{graphicx} 

\usepackage{amssymb, latexsym, amsmath, amscd, array, epigraph,
setspace    
}

\usepackage
[hypertexnames=false]
{hyperref}

\theoremstyle{definition}

\numberwithin{equation}{section}

\newcommand\N {{\mathbb N}} 

\newcommand\R {{\mathbb R}}
\newcommand\Q {{\mathbb Q}}

\newcommand\Los{{\L}o{\'s}}

\newcommand\astr{{}^\ast\hspace*{-0.5pt}\R}

\author[J.B.]{Jacques Bair}\address{J. Bair, HEC-ULG, University of
Liege, 4000 Belgium}\email{j.bair@ulg.ac.be}

\author[P.B.]{Piotr B\l{}aszczyk}\address{P. B\l{}aszczyk, Institute
of Mathematics, Pedagogical University of Cracow,
Poland}\email{pb@up.krakow.pl}

\author[R.E.]{Robert Ely}\address{R. Ely, Department of Mathematics,
University of Idaho, Moscow, ID 83844 US}\email{ely@uidaho.edu}

\author[V.H.]{Val\'erie Henry}\address{V.\;Henry, Department of
Mathematics, University of Namur, 5000
Belgium}\email{vhen@math.fundp.ac.be}

\author[V.K.]{Vladimir Kanovei} \address{V. Kanovei, IPPI, Moscow, and
MIIT, Moscow, Russia} \email{kanovei@rambler.ru}

\author[K.K.]{Karin U. Katz}\address{K. Katz, Department of
Mathematics, Bar Ilan University, Ramat Gan 52900 Israel}
\email{katzmik@math.biu.ac.il}

\author[M.K.]{Mikhail G. Katz}\address{M. Katz, Department of
Mathematics, Bar Ilan University, Ramat Gan 52900 Israel}
\email{katzmik@macs.biu.ac.il}

\author[T.K.]{Taras Kudryk} \address{T. Kudryk, Department of
Mathematics, Lviv National University, Lviv, Ukraine}
\email{kudryk@mail.lviv.ua}

\author[S.K.]{Semen S. Kutateladze}\address{S. Kutateladze, Sobolev
Institute of Mathematics, Novosibirsk State University, Russia}
\email{sskut@math.nsc.ru}

\author[T.M.]{Thomas McGaffey}\address{T. McGaffey, Rice University,
US}\email{thomasmcgaffey@sbcglobal.net}

\author[T. M.]{Thomas Mormann}\address{T. Mormann, Department of Logic
and Philosophy of Science, University of the Basque Country UPV/EHU,
20080 Donostia San Sebastian, Spain} \email{ylxmomot@sf.ehu.es}

\author[D.S.]{David M. Schaps}\address{D. Schaps, Department of
Classical Studies, Bar Ilan University, Ramat Gan 52900 Israel}
\email{dschaps@mail.biu.ac.il}

\author[D.S.]{David Sherry}\address{D. Sherry, Department of
Philosophy, Northern Arizona University, Flagstaff, AZ 86011,
US}\email{david.sherry@nau.edu}

\begin{document}


\thispagestyle{empty}

%
\title{Cauchy, infinitesimals and ghosts of departed
\emph{quantifiers}}

\begin{abstract}
Procedures relying on infinitesimals in Leibniz, Euler and Cauchy have
been interpreted in both a Weierstrassian and Robinson's frameworks.
The latter provides closer proxies for the procedures of the classical
masters.  Thus, Leibniz's distinction between assignable and
inassignable numbers finds a proxy in the distinction between standard
and nonstandard numbers in Robinson's framework, while Leibniz's law
of homogeneity with the implied notion of equality up to negligible
terms finds a mathematical formalisation in terms of standard part.
It is hard to provide parallel formalisations in a Weierstrassian
framework but scholars since Ishiguro have engaged in a quest for
ghosts of departed \emph{quantifiers} to provide a Weierstrassian
account for Leibniz's infinitesimals.  Euler similarly had notions of
equality up to negligible terms, of which he distinguished two types:
geometric and arithmetic.  Euler routinely used product decompositions
into a specific infinite number of factors, and used the binomial
formula with an infinite exponent.  Such procedures have immediate
hyperfinite analogues in Robinson's framework, while in a
Weierstrassian framework they can only be reinterpreted by means of
paraphrases departing significantly from Euler's own presentation.
Cauchy gives lucid definitions of continuity in terms of
infinitesimals that find ready formalisations in Robinson's framework
but scholars working in a Weierstrassian framework bend over backwards
either to claim that Cauchy was vague or to engage in a quest for
ghosts of departed quantifiers in his work.  Cauchy's procedures in
the context of his 1853 sum theorem (for series of continuous
functions) are more readily understood from the viewpoint of
Robinson's framework, where one can exploit tools such as the
pointwise definition of the concept of uniform convergence.  As case
studies, we analyze the approaches of Craig Fraser and Jesper L\"utzen
to Cauchy's contributions to infinitesimal analysis, as well as
Fraser's approach toward Leibniz's theoretical strategy in dealing
with infinitesimals.  The insights by philosophers Ian Hacking and
others into the important roles of contextuality and contingency tend
to undermine Fraser's interpretive framework.  Keywords:
historiography; infinitesimal; Latin model; butterfly model; law of
continuity; ontology; practice; Cauchy; Leibniz.  MSC: 01A45, 01A55,
01A85

\end{abstract}

\maketitle

\tableofcontents

\epigraph{\tiny Some folks see quantifiers when there are none,

Yet others miss infinitesimals when there are
some.

Of history each one is a subversive brat,

Aping the great triumvirate,

While ignoring the insights of Robinson.

(A folk limerick)}

\section{Introduction}

Any given period of the past was characterized by certain practices
and procedures of mathematical researchers.  Such practice may possess
affinities with the practice of modern mathematics, whether of
Weierstrassian or other variety.  Practice reflects something
fundamental about the subject as tackled and experienced by working
mathematicians that is not always expressed in formal foundational
and/or ontological work.  Lurking within past practice are unexplored
potentialities or possibilities.  Some such possibilities did not
develop fully due to contingent factors, but might have developed.
Recent mathematical work can help us in doing history of mathematics
by elucidating potentialities within past mathematical practice.

Our approach differs from that of Barabashev \cite{Ba97} who develops
an \emph{absolute} thesis in favor of \emph{presentism} in the
historiography of mathematics.  We are \emph{not} arguing that one
should use modern mathematics to interpret the past.  While we are not
making any such \emph{absolute} claim, the thesis we wish to develop
is a \emph{relative} one: to the extent that historians are already
using a modern framework to interpret Leibniz, Euler, and Cauchy,
namely a Weierstrassian one, we argue that a better fit for the
\emph{procedures} of these masters is provided by a modern
infinitesimal framework.  Our thesis is \emph{not} that one should
apriori use a modern framework, but that such a modern infinitesimal
framework is preferable to the one traditionally trained historians
often rely upon.

We argue, following the familiar Lakatosian dictum, that history of
mathematics without philosophy of mathematics is blind.  More
specifically, a historiography of mathematics that fails to pay
attention to the dichotomy of procedure versus ontology is flawed.
Following Ian Hacking, we point out that in the historical evolution
of mathematics the amount of contingency is greater than is often
thought.  In particular, this undermines the claim that the historical
development of analysis led to modern classical Weierstrassian
analysis implying that the latter must serve as a \emph{primary
reference point} and interpretative framework.

We compare distinct methodologies in the historiography of mathematics
and point out some pitfalls inherent in neglecting the distinction of
procedure versus ontology, with application to the history of
infinitesimal analysis.  As case studies we analyze Jesper L\"utzen's
approach to Cauchy and Craig Fraser's approach toward both Leibniz's
theoretical strategy in dealing with infinitesimals and Cauchy's
contributions to infinitesimal analysis.  The insights by philosophers
Ian Hacking and others into the important roles of contextuality and
contingency tend to undermine L\"utzen's and Fraser's interpretive
framework.

In 2013, Bair et al.\;\cite{13a} presented some elements of the
history of infinitesimal mathematics often unseen in the received
approach to the history of analysis, based as it is upon a default
Weierstrassian foundation taken as a \emph{primary point of
reference}.  Fraser responded in 2015 with a number of criticisms in
\cite{Fr15}, but his position suffers from some of the very
misconceptions analyzed in \cite{13a}.  These include an insufficient
attention to crucial distinctions such as \emph{practice and
procedure} versus \emph{ontology} as developed by Benacerraf
\cite{Be65}, Quine \cite{Qu68}, and others.  The panoramic nature of
Fraser's survey enables a panoramic overview of the problematic
aspects of the received historiography of mathematical analysis as it
is practiced today.

We spell out some important dichotomies in Section~\ref{sthree}.  We
start with an analysis of L\"utzen's approach in Section~\ref{stwo}
and proceed to Fraser's starting in Section~\ref{srole}.

\section{Dichotomies}
\label{sthree}

Two dichotomies are useful in analyzing the history of analysis:
Archimedean versus Bernoullian, on the one hand, and procedure versus
ontology, on the other.  We will analyze these dichotomies in more
detail in this section.

\subsection{Parallel tracks}
\label{sparallel}

We propose a view of the history of analysis as evolving along
separate, and sometimes competing, tracks.  These are the A-track,
based upon an Archimedean continuum; and B-track, based upon what
could be termed a Bernoullian (i.e., infinitesimal-enriched)
continuum.  Scholars attribute the first systematic use of
infinitesimals as a foundational concept to Johann Bernoulli.  While
Leibniz exploited both infinitesimal methods and \emph{exhaustion}
methods usually interpreted in the context of an Archimedean continuum
(see Bascelli et al.\;\cite{14a}), Bernoulli never wavered from the
infinitesimal methodology.%
\footnote{To note the fact of such systematic use by Bernoulli is not
to say that Bernoulli's foundation is adequate, or that it could
distinguish between manipulations with infinitesimals that produce
only true results and those manipulations that can yield false
results.}

Historians often view the work in analysis from the 17th to the middle
of the 19th century as rooted in a background notion of continuum that
is not punctiform.%
\footnote{Historians use the term \emph{punctiform} to refer to a
continuum thought of as being made out of points.  The term has a
different meaning in modern topology that does not concern us here.}
This necessarily creates a tension with modern, punctiform theories of
the continuum, be it the A-type set-theoretic continuum as developed
by Cantor, Dedekind, Weierstrass, and others; or B-type continua as
developed by Edwin Hewitt \cite{He48} in 1948, Jerzy
\Los{}\;\cite{Lo55} in 1955, Abraham Robinson \cite{Ro66} in 1966, and
others.  How can one escape a trap of presentism in interpreting the
past from the viewpoint of set-theoretic foundations commonly accepted
today, whether of type A or\;B?

\subsection{Procedure versus ontology}
\label{sprocedure}

A possible answer to the query formulated at the end of
Section~\ref{sparallel} resides in a distinction between procedure and
ontology.  In analyzing the work of Fermat, Leibniz, Euler, Cauchy,
and other great mathematicians of the past, one must be careful to
distinguish between
\begin{enumerate}
\item
its syntactic aspects, i.e., actual mathematical practice involving
procedures and inferential moves,  and, 
\item
semantic aspects related to the justification, typically in a
set-theoretic foundational framework, of entities like the points of
the continuum, i.e., issues of the ontology of mathematical entities
like numbers or points.
\end{enumerate}
In his work on Euler and Cauchy, historian Detlef Laugwitz was careful
not to attribute modern set-theoretic constructions to work dating
from before the heroic 1870s, and focused instead on their procedures
(concerning Laugwitz see Section~\ref{s6}).  This theme is developed
further in Section~\ref{s32}.

On the procedure/entity distinction, Abraham Robinson noted that
``from a formalist point of view we may look at our theory
syntactically and may consider that what we have done is to introduce
new \emph{deductive procedures} rather than new mathematical
entities'' \cite[p.\;282]{Ro66} (emphasis in the original).  Kanovei
et al.\;\cite{18d} analyze Robinson's answer to a long-standing
challenge by Felix Klein and Abraham Fraenkel.

Whereas the article on Cauchy's sum theorem by Bascelli et
al.\;\cite{18b} focuses on the details of that result of Cauchy and
argues that Robinson's framework provides the best avenue for
interpreting it, the present text takes a broader lens and identifies
many problems in contemporary historiography of mathematics, using the
work of L\"utzen and Fraser as examples.

\section{``We miss our quantifiers''}
\label{stwo}

In a chapter contributed to a 2003 collection entitled \emph{A History
of Analysis}, L\"utzen writes:
\begin{quote}
(1) To the modern eye Cauchy's definitions may seem wordy, vague and
not particularly rigorous.  (2) We miss our \emph{quantifiers},
our~$\varepsilon$'s,~$\delta$'s and~$N$'s, and in most cases our
inequalities.  \cite[p.\;161]{Lu03} (emphasis and sentence numbering
added)
\end{quote}
Lest a hasty reader should interpret L\"utzen's comment as a frank
acknowledgment of the futility of a quest for ghosts of departed
\emph{quantifiers} in Cauchy, L\"utzen continues:
\begin{quote}
(3) However, as has been pointed out in particular by Grabiner
\ldots{} all these ingredients are clearly present when Cauchy starts
using his concepts in proofs.  (4)\;In particular the complicated
proofs are strikingly \emph{modern} in appearance.  (ibid.) (emphasis
and sentence numbering added)
\end{quote}
The clause ``all these ingredients'' in sentence\;(3) refers to the
$\varepsilon$'s,~$\delta$'s and~$N$'s mentioned in sentence\;(2).
What L\"utzen is claiming is that the~$\varepsilon$'s,~$\delta$'s and
$N$'s appear in Cauchy's proofs, thereby making them ``modern in
appearance'' as per sentence (4).  It is worth pondering the
implications of L\"utzen's apparent assumption that the adjective
\emph{modern} necessarily means~$\varepsilon,\delta$.  Such an
assumption is indicative of a view of analysis as a teleological
process with a specific implied modern outcome, which is decidedly not
that of modern \emph{infinitesimals}.  To escape the trap of
presentism, L\"utzen could have used the expression (in the
terminology of Section~\ref{sparallel}) ``strikingly A-track modern''%
\footnote{In fact L\"utzen appears to interpret the title of the
collection his essay appeared in as referring to ``$A$-History of
Analysis.''}
or some such equivalent turn of phrase.  But being ``strikingly
B-track modern'' is apparently not an option in L\"utzen's book.

Yet Sinkevich pertinently points out that Cauchy's proofs all lack a
characteristic feature of a modern epsilontic proof (``for
each~$\varepsilon$ there exists~$\delta$, etc.''), namely exhibiting
an explicit functional dependence of~$\delta$ on~$\varepsilon$; see
Sinkevich~\cite{Si16}.  Attempts to provide A-type exhaustion
justifications for results in analysis do not originate with Cauchy
and go back at least to Leibniz.  See e.g., Knobloch \cite{Kn02} for
an analysis of an exhaustion procedure in Leibniz, and Bl\aa sj\"o
\cite{Bl17} for an analysis of its limitations.

\subsection{Cauchy's variables}

L\"utzen claims that ``Cauchy writes that a variable having~$0$ as its
limit \emph{becomes} infinitely small'' \cite[p.\;163]{Lu03}.
However, L\"utzen's paraphrase of Cauchy's definition misstates the
order of the significant clauses in Cauchy's own formulation.  In
fact, Cauchy's comment about the \emph{limit} appears not
\emph{before} (as in L\"utzen) but \emph{after} Cauchy's 1823
definition of infinitesimal:
\begin{quote}
Lorsque les valeurs num\'eriques successives d'une m\^eme variable
d\'ecroissent ind\'efiniment, de mani\`ere \`a s'abaisser au-dessous
de tout nombre donn\'e, cette variable devient ce qu'on nomme un
\emph{infiniment petit} ou une quantit\'e infiniment petite.  Une
variable de cette esp\`ece a z\'ero pour limite.%
\footnote{Translation: ``When the successive numerical values [i.e.,
absolute values] of the same variable diminish indefinitely in such a
way as to dip below each given number, this variable becomes what one
calls an \emph{infinitely small} or an infinitely small quantity.  A
variable of this type has limit zero.''}
\cite[p.\;4]{Ca23a} (emphasis in the original)
\end{quote}
This indicates that \emph{limit} is defined in terms of a
\emph{variable quantity} that becomes arbitrarily small, rather than
vice versa.  Indeed Cauchy wrote earlier:
\begin{quote}
Lorsque les valeurs successivement attribu\'ees \`a une m\^eme
variable s'approchent ind\'efiniment d'une valeur fixe, de mani\`ere
\`a finir par en diff\'erer aussi peu que l'on voudra, cette
derni\`ere est appel\'ee la \emph{limite} de toutes les
autres.%
\footnote{Translation: ``When the values successively attributed to
the same variable approach indefinitely a certain fixed value, in such
a way as to differ from it by as little as one wishes, the latter is
called the \emph{limit} of all the others.''}
\cite[p.\;1]{Ca23a} (emphasis in the original)
\end{quote}
The italics on \emph{limite} indicate that this is the concept Cauchy
is in the process of defining.  Thus Cauchy defines both \emph{limits}
and \emph{infinitesimals} in terms of the primitive notion of
\emph{variable quantity}, rather than using \emph{limit} as a
fundamental concept.  L\"utzen goes on to claim that
\begin{quote}
[t]he standard interpretation, also adopted by Grabiner, is that the
limit concept is the central one and infinitesimals only enter as
useful abbreviations for variables having the limit zero
\cite[p.\;163]{Lu03}
\end{quote}
but provides no evidence in Cauchy's texts to back up his ``standard
interpretation'' concerning \emph{limits}.

\subsection{L\"utzen versus Laugwitz--Robinson}
\label{s32b}

L\"utzen tackles what he alleges to be a Laugwitz--Robinson
interpretation of Cauchy's infinitesimals:
\begin{quote}
Laugwitz and Robinson claim that Cauchy's variables not only run
through all values that correspond to our modern real numbers, but
also through infinitesimals as well as sums of real numbers and
infinitesimals.  \ldots{} I find such a revaluation [sic] of Cauchy
interesting because it highlights how historians of mathematics
unconsciously read modern post-Weierstrassian ideas into Cauchy's
work.  \cite[p.\;164]{Lu03}
\end{quote}
L\"utzen claims that Laugwitz and Robinson ``unconsciously read modern
post-Weierstrassian ideas into Cauchy's work.''  To quote a perceptive
adage of the sages of the Talmud, ``\emph{kol haposel, bemumo hu
posel}'' (whoever [seeks to] disqualify [others], it is in his [own]
blemish that he [seeks to] disqualify [them]).  As we will see in
Section~\ref{s34b}, it is actually L\"utzen himself who reads
post-Weierstrassian notions into Cauchy.

Characteristically, L\"utzen does not actually offer any quotes from
either Laugwitz or Robinson to back up his claims.  What L\"utzen
presents as the Laugwitz--Robinson reading is the idea that Cauchy's
variable quantities or magnitudes go not merely through ordinary
values but also infinitesimal values.

L\"utzen goes on to point out that if the variable quantities already
go through infinitesimal values then there would be no need to define
infinitesimals afterwards as Cauchy does, since they are already in
the picture: ``it is hard to explain why infinitesimals are later
defined \emph{afterwards} as variable quantities tending to zero''
\cite[p.\;164]{Lu03} (emphasis in the original) and indeed Cauchy goes
on to define infinitesimals only \emph{after} discussing variable
quantities.  QED for L\"utzen's refutation of Laugwitz--Robinson.

It is not an accident that L\"utzen fails actually to \emph{quote}
either Laugwitz or Robinson on this.  In fact, the interpretation he
seeks to refute is due to neither Laugwitz nor Robinson but
rather\ldots{} Fisher, who wrote in 1978:
\begin{quote}
In the \emph{Pr\'eliminaires} to the \emph{Cours}, Cauchy says: ``When
the successive numerical values of the same variable decrease
indefinitely, in such a way as to fall below any given number, this
variable becomes what one calls \emph{an infinitely small} (\emph{un
infiniment petit}) or an \emph{infinitely small} quantity. A variable
of this kind has zero for limit'' \ldots{} \cite[p.\;316]{Fi78}
\end{quote}
Fisher concludes as follows:
\begin{quote}
The first sentence \emph{could} mean that when the numerical (i.e.,
absolute) values of a variable decrease in such a way as to be less
than any positive number, then the variable takes on infinitesimal
values.  (ibid.)
\end{quote}
L\"utzen may be correct in questioning Fisher's interpretation of
Cauchy's variables, but L\"utzen has not refuted either Laugwitz or
Robinson whom he claims to refute, since such an interpretation is
nowhere to be found in either Laugwitz or Robinson.  L\"utzen merely
introduced further confusion into the subject, and engaged in
regrettable demonisation of Laugwitz and Robinson;
cf.~Section~\ref{s6}.

L\"utzen concludes as follows:
\begin{quote}
If the nonstandard reading of Cauchy is correct, `magnitudes' should
be ordered in a non-Archimedean way and this clashes with the
Euclidean theory.  \cite[p.\;164]{Lu03}
\end{quote}
It is not entirely clear what L\"utzen means by ``Euclidean theory.''
If he means to imply that no non-Archimedean phenomena can be found in
\emph{The Elements}, this is mostly true, but technically incorrect
because we find discussions of horn angles which violate the
Archimedean property when compared to rectilinear angles.  Perhaps by
``Euclidean theory'' L\"utzen means Property\;V.4 of \emph{The
Elements} which is closely related to what we refer today as the
Archimedean property.  If so, L\"utzen seems to claim that Cauchy
would not countenance exploiting infinitesimals which are in blatant
violation of ``Euclidean theory'' a.k.a.\;Property\;V.4.  How
plausible is such a claim?

L\"utzen himself acknowledges that at least Cauchy's \emph{second}
definition of continuity does exploit infinitesimals (see
Section~\ref{s33}) that evidently ``clash[] with the Euclidean
theory'' a.k.a.\;V.4.  Whether magnitudes and quantities go through
infinitesimal values \emph{before} or \emph{after} Cauchy gets around
to defining infinitesimals in terms thereof, there is an irreducible
clash with V.4 (and of course Leibniz already pointed out that
infinitesimals violate V.4; see Bascelli et al.\;\cite{16a}, Bair et
al.\;\cite{17a}).  L\"utzen's allegation that this a shortcoming of
what L\"utzen refers to as a ``nonstandard reading of Cauchy'' is
baseless.

In the Robinson--Laugwitz interpretation, Cauchy starts with variable
quantities which, as far as the 1821 book is concerned, are sequences.
Cauchy says that a sequence getting smaller and smaller \emph{becomes}
an infinitesimal.  Robinson--Laugwitz read this as saying that a
sequence \emph{generates} an infinitesimal which can be expressed in
modern mathematical terminology by saying that infinitesimal~$\alpha$
is an \emph{equivalence class} of sequences, but of course Cauchy had
neither sets nor equivalence classes.  However, the progression of
lemmas and propositions on polynomials and rational functions
in~$\alpha$ that Cauchy gives in the 1821 book does indicate that he
is interested in the asymptotic behavior of the sequence generating an
infinitesimal~$\alpha$; see Borovik--Katz\;\cite{12b} for details.
Thus, a change in finitely many terms of the sequence would not affect
the asymptotic behavior.  A function~$f$ is applied to an
infinitesimal~$\alpha$ by evaluating it term-by-term on the terms of a
sequence generating~$\alpha$.  Thus, if~$\alpha$ is represented by the
sequence~$(1/n)$ then the infinitesimal~$\alpha^2$ will be
incomparably smaller than alpha and will be represented by the
sequence~$(1/n^2)$.

\subsection
{L\"utzen on Cauchy's definitions of continuity}
\label{s33}

L\"utzen goes on to make the following remarkable claim concerning
Cauchy's definitions of continuity:
\begin{quote}
Cauchy actually gives two definitions, first one without
infinitesimals and then one using infinitesimals. \ldots{} The first
definition very clearly specifies a value of the variable~$x$ and
states that~$f(x+\alpha)-f(x)$ tends to zero with~$\alpha$.
\cite[p.\;166]{Lu03}.
\end{quote}
L\"utzen's claim is puzzling all the more so since Cauchy's
definitions are conveniently reproduced in L\"utzen's own text, and
Cauchy clearly describes his~$\alpha$ as ``an infinitely small
increase'' \cite[p.\;159]{Lu03}.  It is therefore difficult, though
not impossible (see below), to maintain as L\"utzen does that Cauchy's
first definition was ``without infinitesimals.''  Cauchy's second
definition is even more straightforward in exploiting an infinitely
small increase~$\alpha$:
\begin{quote}
[Cauchy's second definition] \emph{the function~$f(x)$ will remain
continuous with respect to~$x$ between the given limits if, between
these limits an infinitely small increase [i.e.,~$\alpha$] in the
variable always produces an infinitely small increase in the function
itself}.  \cite[p.\;160]{Lu03} (emphasis in the original)
\end{quote}
The reason we wrote above that it is difficult though not impossible
to maintain as L\"utzen does that Cauchy's \emph{first} definition is
``without infinitesimals'' is because of the possibility of falling
back on a Weierstrassian interpretation of Cauchy's
infinitesimal~$\alpha$ in terms of either sequences or ghosts of
departed quantifiers, as other historians have done.  But if it is
L\"utzen's intention that Cauchy's~$\alpha$ is only a placeholder for
reassuringly Weierstrassian things, then what are we to make of
L\"utzen's acknowledgment on page\;166 (cited above) that Cauchy's
\emph{second} definition \emph{does} exploit infinitesimals?
Apparently the same infinitesimal~$\alpha$ is used in both
definitions, and if~$\alpha$ were only a placeholder in the first
definition, then~$\alpha$ would be a mere placeholder in the second
definition, as well.  Either way you look at it, L\"utzen's reading of
Cauchy's definitions of continuity is problematic.

At any rate we wholeheartedly agree with L\"utzen's assessment that
Cauchy's \emph{second} definition of continuity quoted above does
exploit infinitesimals.  If so, then L\"utzen would apparently have to
agree with us that classical modern analysis in a Weierstrassian
framework is not necessarily the best starting point for understanding
at least some of Cauchy's procedures, and therefore should not
necessarily serve as a primary point of reference when interpreting
Cauchy's second definition of continuity, contrary to Fraser's stated
position; see Sections~\ref{srole} through \ref{s6}.

\subsection{Did Cauchy need Dedekind?}
\label{s34b}

In this Section we examine the following startlingly ahistorical claim
by L\"utzen:
\begin{quote} 
In modern treatment [the convergence of a Cauchy sequence] is derived
from the completeness of the real numbers (or is even taken as the
definition of completeness) which must be either postulated as an
axiom or obtained from a construction of the real numbers.  Cauchy's
work does not contain either way out and he could not have appealed to
the underlying concept of magnitude because Euclid does not have
axioms which ensure the completeness of his magnitudes.

This missing account of completeness is a \emph{fundamental lacuna}
which appears in several other places in Cauchy's analysis, in
particular in his proof of the intermediate value theorem\ldots{}
\cite[p.\;167-168]{Lu03} (emphasis added)
\end{quote}
According to L\"utzen, it is a ``fundamental lacuna'' of Cauchy's
proof of intermediate value theorem that the result relies on
completeness and completeness could not have been provided by Cauchy.
Why is that?  Says L\"utzen, because there are two ways of doing that:
either axiomatic or construction.  L\"utzen claims that Cauchy did not
have a construction, and as far as axioms go, Euclid did not have
anything on completeness.

However, L\"utzen'a analysis is misleading not merely because Simon
Stevin already gave us unending decimals centuries earlier, but also
because Laugwitz explicitly points out in his papers (see e.g.,
\cite[p.\;202]{La89}) that Cauchy did not need a construction of the
reals because he had unending decimal expansions.

Cauchy's argument for finding a zero of a continuous function that
changes sign works because it gives a procedure for finding an
unending decimal expansion of such a zero, without hifalutin'
embellishments concerning completeness, which is certainly a useful
notion but was only introduced by Dedekind half a century later.

Thus, in the context of unending decimals Cauchy \emph{did} have a
suitable construction (contrary to L\"utzen's claim concerning an
alleged \emph{fundamental lacuna}) especially if one uses recursive
subdivision into~$10$ equal parts as Stevin had done, which would
directly produce an unending decimal expansion of a required zero of
the function; see B\l aszczyk et al.\;\cite[Section\;2.3]{13b} for
details.

The Whig history aspect of L\"utzen's analysis is ironic since this is
precisely the shortcoming he seeks to pin on Robinson and Laugwitz;
see Section~\ref{s32b}.

\emph{Without proposing any textual analysis} of Cauchy's proof,
L\"utzen concludes, apparently on the grounds of historical hindsight
deduced from Dedekind backwards, that Cauchy's proof \emph{could not
possibly be correct}, because Dedekind hand't clarified the relevant
mathematical notion of completeness yet.  But, as we saw, Cauchy did
not need Dedekind simply because the assertion (or tacit assumption)
that each unending decimal results in a legitimate real number (1) is
equivalent to the completeness, and furthermore (2) could hardly be
questioned by any sound pre-Weierstrassian analyst.  The axiomatic
formalisation of completeness was not a dramatic new discovery but
rather a workhorse systematisation of common views in a new
revolutionary language.

\section{Role of modern theories}
\label{srole}

The relevance of modern theories to interpreting the mathematics of
the past is a knotty issue that has caused much ink to be spilled.
Without attempting to resolve it, we offer the following thoughts to
help clarify the issue.

\subsection
{On the relevance of theories for history of mathematics}
\label{s13}

Let us begin by considering, as a thought experiment, the following
passage.
\begin{quote}
(B) The relevance of modern \emph{Archimedean} theories to an
historical appreciation of the early calculus is a moot point.  It is
doubtful if it is possible or advisable to reconstruct a past
mathematical subject in a way that conforms to modern theories and is
also consistent with how a practitioner of the period would have
worked.  Such a project runs the obvious risk of imposing one's own
interests and conceptions on the past subject.  It is possible that
such an endeavor will end up with something of intellectual interest
and mathematical value, but it is unlikely that it will constitute a
\emph{significant contribution to history}.
\end{quote}
The source of the passage (B) will be identified below.  For the time
being we note the following.  The passage amounts to a sweeping
dismissal of the modern, punctiform, Archimedean framework (as
developed by Cantor, Dedekind, Weierstrass and others) as a primary
point of reference for historical work dealing with prior centuries.

A reader might react with puzzlement to such a sweeping dismissal,
feeling that the grounds for it are too generic to be convincing.
Granted the punctiform nature of the modern continuum is different
from the historical views of the continuum as found for example in
Leibniz, Euler, and Cauchy.%
\footnote{For a study of Stevin, see Katz--Katz \cite{12c} and
B\l{}aszczyk et al.\;\cite[Section\;2]{13b}; for Fermat, see Katz et
al.\;\cite{13e} and Bair et al.\;\cite{18g}; for Gregory, see Bascelli
et al.\;\cite{17b}.}
However, the reader may feel that the passage above ends up throwing
out the \emph{baby} with the (punctiform) \emph{bathwater}.

To avoid such an overreaction one may argue as follows.  For all the
set-theoretic novelty of the modern view of the Archimedean continuum,
surely important insight is to be gained from modern axiomatisations
of the \emph{procedures} in Euclid's \emph{Elements}, resulting in a
greater appreciation, and understanding, of this classical work (on
Euclid see further in Section~\ref{s4}).

Similarly, historians have often endeavored to show how the lucidity
of the modern Archimedean \emph{limit} concept serves as a benchmark
that allows us to appreciate the progress being made over the decades
and centuries in clarifying the nature of the \emph{procedures}
involving this concept that were employed in the 17th through the 19th
centuries in solving problems of analysis.

To save the \emph{baby}, the reader may feel that a distinction needs
to be made between, on the one hand, the modern set-theoretic
justification of the entities involved, such as that of the complete
Archimedean ordered field; and, on the other, the \emph{procedures}
and inferential moves exploited in mathematical arguments, as already
discussed in Section~\ref{sprocedure}.  The former (set-theoretic
justification) pertains to the domain of the ontology of mathematical
entities that apparently had little historical counterpart prior to
1870.  Meanwhile, the latter (\emph{procedures}) are, on the contrary,
indispensable tools that allow us to clarify and appreciate the
inferential moves as found in the work of the mathematicians of the
past.

The distinction between procedure and ontology is a tool implicitly
used in \cite{Fr15} (see Section~\ref{s16}).  Similarly, Ferraro and
Panza, both shortlisted in \cite[p.\;27]{Fr15} among the authorities
on the 18th century, declare a \emph{failure} of Lagrange's program of
founding analysis on power series in the following terms:
\begin{quote}
[Lagrange's books] are part of a foundationalist agenda. The fact that
this agenda was never really accepted by Lagrange's contemporaries
contrasts with another fact: that it was the most careful attempt to
integrate the calculus within the program of 18th-century algebraic
analysis.  Its \emph{failure} is then also the \emph{failure} of this
ambitious program. \cite[p.\;189, Conclusions]{FP} (emphasis added)
\end{quote}
One might ask, in what sense did Lagrange's program \emph{fail}, if
not in the modern sense of being unable to account for non-analytic
functions?  When restricted to analytic ones, Lagrange's program
indeed succeeds.  Ferraro and Panza focus on Lagrange's \emph{purity
of the method} but ultimately use a modern yardstick to offer an
evaluation (namely, an alleged \emph{failure}) of Lagrange's program,
without necessarily falling into a presentist trap of attributing a
punctiform continuum to Lagrange.

We hasten to add that the passage (B) cited at the beginning of the
current Section~\ref{s13} is Fraser's only up to a single, but very
significant, modification.  We changed his adjective
\emph{non-Archimedean} to \emph{Archimedean}.  Fraser's original
passage reads as follows.

\subsection{Fraser's original passage}
\label{s42}
(A) ``The relevance of modern non-Archimedean theories to an
historical appreciation of the early calculus is a moot point.  It is
doubtful if it is possible or advisable to reconstruct a past
mathematical subject in a way that conforms to modern theories and is
also consistent with how a practitioner of the period would have
worked.  Such a project runs the obvious risk of imposing one's own
interests and conceptions on the past subject.  It is possible that
such an endeavor will end up with something of intellectual interest
and mathematical value, but it is unlikely that it will constitute a
significant contribution to history.''  \cite[p.\;43]{Fr15}.

\subsection{Analysis of Fraser's original}
\label{s32}

Fraser's original passage (A) appearing in Section~\ref{s42} amounts
to a sweeping dismissal of the modern non-Archimedean framework as
developed by Robinson, Laugwitz, and others as a primary point of
historical reference.

Again, the reader might react with puzzlement to such a sweeping
dismissal on what appear to be excessively generic grounds.  The
reader might object to Fraser's throwing out the baby with the
(set-theoretic) bathwater when he ignores the distinction between
modern foundational \emph{justification} and classical
\emph{procedures}.  

The upshot of our thought experiment should be clear by now.  For
Fraser, neither the modern Archi\-me\-dean nor the modern
non-Archimedean theories should be relevant to a historical
appreciation of the early calculus.  Rather, for him, history of
mathematics stands alone and aloof.  As a result it risks becoming,
one may say, a vain (in both senses) exercise of ``\emph{l'Histoire,
c'est Moi}'' type; see Katz \cite{Ka15}.

In his passage (A), Fraser implicitly argues in favor of a
\emph{splendid isolation} for a historiography of mathematics;
according to him, it is to be pursued independently of any influence
coming from modern mathematical theories whatsoever, whether of the
Archimedean or the Bernoullian variety.  This naturally extends to a
claim that history of mathematics should be pursued independently of
any theoretical modern influences whatsoever: after all, modern
mathematics comes along with, like it or not, a lot of
\emph{ideological baggage} concerning problems of what mathematics is,
how it is related to other sciences, etc.  

Yet the contention that one's own conception of the history of
mathematics would be free of any such baggage is hardly convincing
(see especially Fraser's remarks on classical analysis as a ``primary
point of reference'' analyzed in Section~\ref{s36}).  Thus, replacing
the somewhat unfriendly expression \emph{ideological baggage} by the
more neutral term \emph{philosophy} we may assert that Fraser's
intended isolationism amounts to a \emph{history of mathematics
without philosophy of mathematics}.

While we cautiously subscribe to Lakatos's dictum concerning the
dependence of history of science on philosophy of science,%
\footnote{Lakatos's insight as formulated in \cite{La71} was grounded
upon earlier historical-philosophical studies such as Hesse
\cite{He62}.}
we don't subscribe to the way in which he applied it to the history of
the early calculus, preferring Laugwitz's analysis as developed in the
pages of \emph{Historia Mathematica} \cite{La87}, \emph{Archive for
History of Exact Sciences} \cite{La89}, and elsewhere.

As philosopher Marx Wartofsky pointed out in his programmatic
contribution \emph{The Relation between Philosophy of Science and
History of Science} \cite{Wa76}, there are many distinct possible
relations between philosophy of science and history of science, some
``more agreeable'' and fruitful than others (ibid., p.\;719ff).  We
need not go into the details of Wartofsky's typology of possible
relations between the two disciplines.  It will suffice to point out
that according to him a fruitful relation between history and
philosophy of science requires a rich and complex \emph{ontology} of
that science.

In the case of mathematics, this means that a fruitful relation
between history and philosophy must go beyond offering an ontology of
the domain over which a certain piece of mathematics ranges (say,
numbers, functions, sets, infinitesimals, structures, or whatever).
Namely, it must develop the ontology of mathematics \emph{as a
scientific theory} itself (ibid., p.\;723).

In the present article, we make a step in this direction by
distinguishing between the (historically relative) \emph{ontology} of
the mathematical objects in a certain historical setting, and its
\emph{procedures}, particularly emphasizing the different roles these
components play in the history of mathematics.  More precisely, our
\emph{procedures} serve as a representative of what Wartofsky called
the \emph{praxis} characteristic of the mathematics of a certain time
period, and our \emph{ontology} takes care of the mathematical objects
recognized at that time.

The urgency of drawing attention \emph{away} from mathematical
foundations (in a narrow ontological sense) and \emph{toward} practice
was expressed by Lawvere:
\begin{quote}
In my own education I was fortunate to have two teachers who used the
term `foundations' in a common-sense way (rather than in the
speculative way of the Bolzano--Frege--Peano--Russell tradition).
This way is exemplified by their work in \emph{Foundations of
Algebraic Topology}, published in 1952 by Eilenberg (with Steenrod),
and \emph{The Mechanical Foundations of Elasticity and Fluid
Mechanics}, published in the same year by Truesdell. The orientation
of these works seemed to be `concentrate the essence of practice and
in turn use the result to guide practice.'  \cite[p.\;213]{La03}
\end{quote}
Related educational issues were analyzed by Katz--Polev \cite{17h}.

Our thought experiment in Section~\ref{s13} comparing the passages (B)
and~(A) helps highlight the problem with Fraser's argument.  Fraser
attempts to argue against using modern Bernoullian frameworks (see
Section~\ref{sparallel}) to interpret the mathematics of the past, but
fails to take into account the crucial distinction between procedure
and ontology.  Therefore his argument in favor of dismissing modern
frameworks is so general as to apply to his own work, and indeed to
much valuable work on understanding the mathematics of the past, and
has the effect of so to speak throwing out the (procedural) baby with
the (ontological) bathwater.  He wishes to reject applications of
Bernoullian frameworks to procedures of the mathematics of the past
while ignoring much of the evidence, and on such generic grounds as to
make his critique untenable.

\subsection{A conceptual gulf}
\label{s34}

A failure to keep in mind crucial distinctions such as that between
procedure and ontology undermines Fraser's evaluation of Robinson's
historical work.  Writes Fraser:
\begin{quote}
The transition from algebraic analysis of the eighteenth century to
Cauchy--Weierstrass analysis of the nineteenth century marked a much
greater discontinuity than did the emergence of nonstandard analysis
out of classical analysis in the second half of the twentieth century.
\cite[p.\;42]{Fr15}
\end{quote}
This is a valid remark (particularly if one omits the mention of
Cauchy) on the nature of post-Weierstrassian analysis, and one
consistent with our comments on the ontology of punctiform continua in
Section~\ref{sthree}.  Fraser continues:
\begin{quote}
Nonstandard analysis is an offshoot of modern analysis and sits
solidly on the modern side of the conceptual gulf opened up by the
\emph{Cauchy--Weierstrass foundation}. (ibid.)
\end{quote}
Fraser again makes a valid point (again modulo Cauchy's role)
concerning the ontology underwriting mathematical entities since
Weierstrass.  But then Fraser makes the following leap:
\begin{quote}
In this respect Robinson and to a lesser degree Lakatos were mistaken
in their assessment of Cauchy.  (ibid.)
\end{quote}
This remark of Fraser's is a \emph{non-sequitur} exacerbated by
Fraser's failure to make explicit the alleged mistakes of Robinson and
Lakatos.  The following comment by Robinson already quoted in
Section~\ref{sprocedure} may help set the record straight:
\begin{quote}
\ldots from a formalist point of view we may look at our theory
syntactically and may consider that what we have done is to introduce
\emph{new deductive procedures} rather than new mathematical entities.
\cite[p.\;282]{Ro66} (emphasis in the original)
\end{quote}
Nor does Fraser provide any detailed examination of Robinson's
analysis of Cauchyan texts (though a brief summary appears in
\cite[p.\;26]{Fr15}).  In \emph{what respect} exactly is Robinson
allegedly ``mistaken in [his] assessment of Cauchy'' as Fraser claims?

Following on the heels of Fraser's comments on the \emph{conceptual
gulf} wrought by set-theoretic foundations that emerged following the
work of Weierstrass, Fraser's remark derives its impetus from the
modern nature of Robinson's framework based as it is on punctiform
models of continua, contrary to pre-1870 work which was generally not
based on such ontology.  Yet Fraser fails to take into account the
fact that in discussing Cauchy, Robinson is talking \emph{procedure},
not \emph{ontology}.%
\footnote{It is worth mentioning that, among other sources, Robinson's
work grew from a reflection upon Skolem's nonstandard models of
arithmetic.  These were developed in \cite{Sk33}, \cite{Sk34},
\cite{Sk55}.  Skolem represents his nonstandard numbers by definable
sequences of integers (the standard numbers being represented by
constant sequences).  A sequence that tends to infinity generates an
infinite number.  This procedure bears analogy to Cauchy's
representation of the B-continuum.  Indeed, Cauchy gives an example of
a variable quantity as a \emph{sequence} at the start of his
\emph{Cours d'Analyse}\;\cite{Ca21}.}

\subsection{Who took Cauchy out of his historical milieu?}
\label{s35}

Having summarily dismissed Robinson, Fraser goes on to yet another
non-sequitur:
\begin{quote}
Rather than try to understand Cauchy as someone who developed within a
given intellectual and historical milieu, they approach the history
from an essentially artificial point of view.  \cite[p.\;42]{Fr15}
\end{quote}
Fraser provides no evidence for his claim that Robinson seeks to take
Cauchy out of his historical milieu.  One may well wonder whether, on
the contrary, it is Fraser who seeks to take Cauchy out of his
historical milieu.  Indeed, Robinson takes Cauchy's infinitely small
at (procedural) face value based on an assumption that Cauchy
understood the term in a sense common among his contemporaries like
Abel and Poisson, as well as a majority of his colleagues at the Ecole
Polytechnique, both mathematicians and physicists, whose work was a
natural habitat for infinitesimals, as documented by
Gilain\;\cite{Gi89}.

Meanwhile, by postulating a so-called \emph{Cauchy--Weierstrass}
foundation (see Section~\ref{s34}), Fraser precisely yanks Cauchy
right out of his historical milieu, and inserts him in the heroic
1870s alongside C. Boyer's great triumvirate.%
\footnote{\label{f1}Historian Carl Boyer described Cantor, Dedekind,
and Weierstrass as \emph{the great triumvirate} in \cite[p.\;298]{Boy}.
The term serves as a humorous characterisation of both A-track
scholars and their objects of adulation.}

Robinson was one of the first to express the sentiment that the
A-framework is inadequate to account for Cauchy's infinitesimal
mathematics.  Grattan-Guinness points out that Cauchy's proof of the
sum theorem is difficult to interpret in an Archimedean framework, due
to Cauchy's use of infinitesimals.  Thus, Grattan-Guinness wrote:
\begin{quote}
This remark [of Cauchy's] is \emph{difficult to interpret} against
[i.e., in the context of] the classification of modes of uniform
convergence given here \ldots{} since~$\alpha$ is an infinitesimally
small increment of~$x$. \cite[p.\;122]{GG70} (emphasis added)
\end{quote}
Recent studies by Nakane \cite{Na14} and S\o{}rensen \cite{So05} have
emphasized the difference between Cauchyan and Weierstrassian notions
of limit.

Meanwhile, a B-track framework enables a better understanding of
Cauchy's procedures in his 1853 text on the convergence of series of
functions \cite{Ca53} (see Bascelli et al.\;\cite{18b} for a detailed
analysis) and other texts where infinitesimals are used in an
essential fashion, such as Cauchy's seminal 1832 text on geometric
probability \cite{Ca32} and his seminal 1823 text on the theory of
elasticity \cite{Ca23b}.  Fraser's misjudgment of Cauchy's historical
role is a predictable consequence of his dogmatic endorsement of
A-track analysis as the \emph{primary point of reference} for
understanding the mathematics of the past, as we discuss in
Section~\ref{s36}.

\subsection{The apple of discord}
\label{s36}

Fraser seeks to distance himself from Boyer's view of the history of
analysis as inevitable progress toward arithmetisation:
\begin{quote}
Since the 1960s there has been a new wave of writing about the history
of eighteenth-century mathematics [that has] charted the development
of calculus without interpreting this development as a first stage in
the inevitable evolution of an arithmetic foundation.  (Fraser
\cite[p.\;27]{Fr15})
\end{quote}
However, the following passage from Fraser's survey reveals the nature
of his misconception:
\begin{quote}
Of course, classical analysis developed out of the older subject and
it remains a \emph{primary point of reference} for understanding the
eighteenth-century theories.  By contrast, nonstandard analysis and
other non-Archimedean versions of calculus emerged only fairly
recently in somewhat \emph{abstruse} mathematical settings that bear
\emph{little connection} to the historical developments one and a
half, two or three centuries earlier.  (ibid.; emphasis added)
\end{quote}
For all his attempts to distance himself from Boyer's idolisation of
the triumvirate,%
\footnote{See note~\ref{f1} on page~\pageref{f1} for an explanation of
the term.}
Fraser here commits himself to a position similar to Boyer's.  Namely,
Fraser claims that modern punctiform A-track analysis is, ``of course,
[!] a primary point of reference'' for understanding the analysis of
the past.  His sentiment that modern punctiform B-track analysis
\emph{bears little connection to the historical developments} ignores
the procedure versus ontology dichotomy.

Arguably, modern infinitesimal analysis provides better proxies for
the procedures of the calculus of the founders, based as it was on a
fundamental Leibnizian distinction between assignable and inassignable
quantities%
\footnote{The distinction goes back to the work of Nicholas of Cusa
(1401--1464), which also inspired Galilio's distinction between
\emph{quanta} and \emph{non-quanta} according to
Knobloch\;\cite{Kn99}.}
underpinning the infinitesimal analysis of the 17th and 18th
centuries, as analyzed in B\l aszczyk et al.\;\cite{17d}.  Meanwhile
scholars since Ishiguro \cite[Chapter\;5]{Is90} have been engaged in a
syncategorematic quest for ghosts of departed quantifiers in Leibniz.

A sentiment of the inevitability of classical analysis is expressed by
Fraser who feels that ``classical analysis developed out of the older
subject and it remains a primary point of reference for understanding
the eighteenth-century theories'' (ibid.); yet his very formulation
involves circular reasoning.  It is only if one takes classical
analysis as \emph{a primary point of reference} that it becomes
plausible to conclude that it \emph{developed out of the older
subject} (and therefore should serve as a primary point of reference).

To comment on Ian Hacking's distinction between the \emph{butterfly
model} and the \emph{Latin model}, we note the contrast between a
model of a deterministic biological development of animals such as
butterflies, as opposed to a model of a contingent historical
evolution of languages such as Latin.

If one allows, with Ian Hacking \cite[pp.\;72--75]{Ha14}, for the
possibility of alternative courses for the development of analysis (a
\emph{Latin model} as opposed to a \emph{butterfly model}), Fraser's
assumption that the historical development necessarily led to modern
classical analysis (as formalized by Weierstrass and others) remains
an unsupported hypothesis.

The butterfly model assumes a teleological or culminationist view of
mathematical development, in which the texts of historical
mathematicians are interpreted in terms of their location on the track
that inevitably culminates in our current dominant formal
understanding of analysis. In contrast, the Latin model treats some
aspects of mathematical development as contingent rather than
inevitable.

By distinguishing A-track from B-track development, we are not
therefore proposing a butterfly model, just with two possible
butterflies at the end.  Rather, we are noting that it is an illusion
that one can interpret historical mathematics without the
interpretation being informed by philosophical commitments or
mathematical frameworks. Fraser and others we have referenced in this
paper have not applied a framework-free interpretation of the works of
mathematicians such as Leibniz, Euler, and Cauchy. They have,
consciously or unconsciously, employed an A-track or Weierstrassian
framework. By focusing attention on mathematical procedure rather than
ontology, we can see that a B-track framework reveals aspects of these
great mathematicians' work that A-track analyses distort or dismiss.
The B-track analysis also provides an alternative to the
claustrophobic presumption that mathematical development is
pre-scripted or determined.

We will respond to Fraser's \emph{little connection} claim in more
detail in Section~\ref{s4}.

\section{Euclid}
\label{s4}

In this section we will respond more fully to Fraser's passage quoted
in Section~\ref{s36}.

\subsection{History and `little connection'} 

Since 1870, both non-Archime\-dean analysis and non-Archimedean
geometry have developed in parallel to classical analysis.  Here one
could mention the work of Otto Stolz and Paul du Bois-Reymond in
analysis, and that of Giuseppe Veronese and David Hilbert in geometry;
see Ehrlich \cite{Eh06}.

Hermann Hankel developed the first modern interpretation of Book~V of
Euclid's \emph{Elements} in his 1876 work \cite{Ha76}.  In this way he
initiated mathematical investigations of the notion of magnitude
carried out further by Stolz, Du Bois-Reymond, H. Weber, and
O. H\"{o}lder; for details see B\l aszczyk \cite{Bl13}.

This line of investigation culminated in 1899/1903 with the first and
second editions of \emph{Grundlagen der Geometrie}, when Hilbert
pioneered an axiom system for an ordered field; see Hilbert
\cite{Hi99}, \cite{Hi03}.  With his axiomatic method, Hilbert
initiated a new type of investigation of a foundational type in
geometry.  However, Hilbert also initiated an axiomatic study of the
real numbers, which is a lesser-known aspect of his work.  Such new
investigations focus, on the one hand, on axiomatic characterisations
of mathematical structures, and on the other, on the descriptive power
of formal languages.

As regards the real numbers, foundational studies were carried out,
among others, by E.\;V.\;Huntington and O.\;Veblen, at the beginning
of the 20th century.  In 1926 Artin and Schreier \cite{AS} and in 1931
Tarski \cite{Ta31} developed the theory of real-closed fields.
Robinson explicitly referred to the ``theory of formally-real fields
of Artin and Schreier'' in \cite[p.\;278]{Ro66}.  Robinson's framework
for infinitesimal analysis inscribes in the tradition of foundational
studies of the real numbers.

Hilbert's foundational studies in geometry are well-known.  While his
\emph{Grundlagen der Geometrie} is primarily a mathematical book, it
is referred to both by mathematicians studying the foundations of
Euclidean geometry and by historians interpreting Euclid's geometry;
see e.g., Heath \cite{He08}, Mueller \cite{Mu81}.  Historians of Greek
mathematics often rely on modern interpretive techniques.  Thus Netz
\cite{Ne99} refers to cognitive science and Kuhn's philosophy.
Naturally also refer to the achievements of modern mathematics.

Now, a given historian's competence and the limitations of the
mathematical techniques he learned, tend to determine what portion of
contemporary mathematics he is able to apply in his historical
studies.%
\footnote{Such limitations are presumably what prompts a historian to
describe a given piece of mathematics as \emph{abstruse}; see
Section~\ref{s36} and Section~\ref{s10}.}
This applies even to the history of Greek mathematics, where the
distance between the source text we refer to and mathematical
techniques we apply in our interpretation is measured in thousands of
years, as compared to hundreds in the case of infinitesimal calculus.
In the next two sections we will provide some examples.

\subsection{Applying modern mathematics to interpreting Euclid} 

The very first proposition of Euclid's \emph{Elements} is
controversial as it relies on the assumption that two circles
occurring in its proof must meet.  Almost every commentator from
Proclus to Mueller comments on this assumption.  Already in the 1645
Claude Richard \cite{Ri45} realized that an additional hypothesis is
required here; see de Risi \cite[section~1.2]{De15}.

\subsubsection{Hilbert's axioms}

Heath claimed that one ``must invoke'' the continuity axiom to fill
this logical gap in Euclid's proof \cite[vol.\;1, p.\;242]{He08}.
However, interpretaters of Euclid's \emph{Elements} following
Hilbert's \emph{Grundlagen der Geometrie} showed that in place of the
continuity axiom, it is sufficient to invoke the so-called Pasch axiom
and circle-circle intersection property, and then Euclid's geometry,
as developed in Books 1--4, can be developed in a plane without the
continuity axiom; see Hartshorne \cite[p.\;110]{Ha00}.

Furthermore, the Pasch axiom can be derived from the fact that the
field of segments is a Pythagorean field, i.e., a field closed under
the operation~$\sqrt{1+a^2}$; see \cite[p.\;187]{Ha00}.

\subsubsection{Modern axiomatic systems}

Modern axiomatic systems that followed, do not appeal explicitly to
the notion of a circle.  Instead, they encode the circle-circle
intersection property in some axioms which apply relations of
congruence, betweenness, and equal distance.  Thus, Karol Borsuk and
Wanda Szmielew modified Hilbert axiom system, and showed that Euclid's
geometry can be developed based on axioms of incidence, congruence,
equal distance and order, including the Pasch axiom, without the
continuity axiom; see \cite[\S\,93]{BS}.

\subsubsection{Tarski's axioms}

In 1959 Alfred Tarski \cite{Ta59} obtained further results using a
system of axioms he developed.  This system can be formulated in first
order logic and enables one to \emph{translate} some metageometrical
results into \emph{ordered field} terms.  Thus, it turns out that
models of Euclidean geometry coincide (up to isomorphism) with
Cartesian planes over a \emph{Euclidean field}, i.e., an ordered field
closed under the square root operation.

In 1974 Szmielew \cite{Sz74} showed that models of geometry satisfying
the Pasch axiom coincide (up to isomorphism) with a Cartesian plane
over a Pythagorean field.  The so called \emph{two-circles} axiom%
\footnote{The \emph{circle} axiom asserts that
$B(abc)\Rightarrow(\exists c')[B(pbc')\wedge D(ac'ac)]$.  The
\emph{two circles} axiom asserts that
\[
 M(ac'b)\wedge M(ab'c)\wedge M(ba'c)\Rightarrow (\exists
q)[L(pc'q)\wedge (D(qb'ab')\vee D(qa'ba'))],
\]
where~$M(abc)$ means that~$b$ is the midpoint, i.e., lies between~$a$
and~$c$, and~$D(abbc)$ while~$L$ is the \emph{colinearity}
relation,~$B$ is the \emph{betweenness} relation, and~$D$ is the
\emph{equal distance} relation.}
provides the geometric analog for the Euclidean field property.
Furthermore,
\begin{enumerate}
\item
the Pasch axiom is a consequence of the \emph{circle} axiom;
\item
the Pasch axiom does not follow from the continuity axiom, and
\item
the circle axiom indeed follows from the continuity axiom;
\end{enumerate}
see respectively Szmielew \cite{Sz70a}, Szmielew \cite{Sz70b}, and
Szczerba \cite{Sz70}.  With regard to this somewhat paradoxical
situation, Szmielew commented as follows: ``At the first glance the
situation seem to be paradoxical, since it is common to say that C
[the circle axiom] is a particular consequence of Co [the continuity
axiom].  In fact, C can not be proved without the use of Co, however,
in the proof of C besides Co also some other axioms of\;E [Euclidean
geometry] are involved'' \cite[p.\;751]{Sz70b}.

Relating modern axioms, on the one hand, to Euclid's axioms,
postulates and theorems, on the other, constitutes a methodological
challenge since that there is no explicit counterpart of betweenness
relation in Euclid's geometry.  Nonetheless, such metageometrical
results arguably contribute both to our interpretation and
appreciation of Euclid's \emph{Elements} and to the book's historical
perception, by enabling new interpretive perspectives, be it
\emph{logical gaps} in Euclid's development or the \emph{role of
diagrams} in his proofs.  This is contrary to Fraser's contention that
analysis of historical texts that relies on modern mathematics cannot
produce a \emph{significant contribution to history}.  While this
methodology is widely accepted, it differs from the one we apply only
in the scope of modern mathematics being employed.

\subsection
{Applying modern mathematics to the real numbers}
 
Our next example concerns the real numbers.  In light of the
foundational studies mentioned above, we have today not only different
versions of the continuity axiom but also some independence results.
We know, for example, that continuity of real numbers can be
characterized by the Dedekind axiom (i.e., no Dedekind cut gives a
gap) or equivalently by the conjunction of \emph{two} conditions:
Cauchy completeness (CC) and the Archimedean axiom (AA).  It is known
that axioms CC and AA are independent, since both the hyperreals and
Levi-Civita fields are non-Archimedean but Cauchy-complete.

With this knowledge we now consider the historical texts. In 1872
neither Cantor, Dedekind, Heine, or M\'eray were aware of AA.\,
However, while Dedekind characterized the real numbers by his axiom
(and from the current perspective it is a correct characterisation),
Cantor, as well as M\'eray and Heine, characterized real numbers by
means of CC alone.  Can we then claim that they (i.e., Cantor, Heine,
and M\'eray) really knew what real numbers are?  Apparently not.

Furthermore, while Cauchy sequences, or \emph{Fundamentalreihen}, were
widely in use at the time, Dedekind's axiom by no means bears any
connection to historical developments that preceded \emph{Stetigkeit
und irrationale Zahlen}.
 
To give a further example, the Intermediate Value Theorem (IVT) is
another equivalent form of the the Dedekind Cut Axiom; see Teismann
\cite[section\;3]{Te13}, and B\l aszczyk \cite[note\;5]{Bl15}.  This
relates to the famous Bolzano proof of the IVT \cite{Bo04}.  Bolzano
had no access to axiomatic characterisations of the real numbers.
Should this modern result be employed in analyzing Bolzano's pamphlet?
We argue that the answer is affirmative, as this result suggests a new
interpretive perspective of seeking another version of the continuity
axiom that he implicitly used in his proof.

As far as Robinson's framework is concerned, we have shown that, from
a historical perspective, it belongs squarely in the tradition of
foundational studies of the real numbers.  However, Robinson's
framework also sheds new light on the relation of \emph{being
infinitely close}, as well as on infinite (or in modern terminology,
\emph{hyperfinite}) sums and products (see Section~\ref{s54}).
Therefore Robinson's framework can be usefully applied in
historiography on a par with Hilbert's axiomatic method and the
results it provides, beyond the strictures of a Weierstrassian
framework useful though it may be.

\section{From general to the specific}

In addition to general claims concerning an alleged \emph{a priori}
inapplicability of Bernoullian frameworks to historical research as
analyzed in Section~\ref{srole}, Fraser's article contains some
specific claims concerning Robinson's framework.  We will analyze one
such claim in this section.

\subsection{Algorithms}
\label{s16}

As we mentioned in Section~\ref{srole}, Fraser himself implicitly
relies on the dichotomy of procedure versus ontology.  Thus, he
writes:
\begin{quote}
Two of the most prominent features of the early calculus -- the
tension between analytic and geometric modes of representation and the
\emph{central place occupied by the algorithm} -- are not reproduced
at all in nonstandard analysis as defining characteristics of the
subject.  \cite[p.\;39]{Fr15} (emphasis added)
\end{quote}
We will examine the ``tensions'' mentioned by Fraser on another
occasion.  What interests us here is Fraser's claim concerning ``the
central place occupied by the algorithm.''  Now \emph{algorithmic}
aspects clearly fall under the rubric of \emph{procedures} as opposed
to \emph{ontology} (see Section~\ref{sprocedure}).  Focus on
algorithms certainly characterized 18th century mathematics.  We will
examine Fraser's claim that they allegedly ``are not reproduced at all
in nonstandard analysis'' in Section~\ref{s52}.

\subsection
{Fraser, algebraic analysis, and the details}

Fraser argues that infinitesimals and the familiar foundational issues
surrounding them were not the primary motivators of the evolution of
analysis.  He further claims that the bulk of 19th-century work
responded to the context of the algebraic-analytic approach
exemplified by Lagrange, in which infinitesimals already played a
secondary role:
\begin{quote}
The question of the logical status of infinitesimals is of
\emph{secondary interest}, from either the perspective of a researcher
in the early eighteenth century or an observer today.  The second
topic concerns the decisive shift to \emph{algebraic analysis} that
occurred in the writings of such figures as Euler and Lagrange in the
second half of the century.  \cite[p.\;28]{Fr15} (emphasis added)
\end{quote}
Fraser continues:
\begin{quote}
The French philosopher Auguste Comte \ldots{} revered Lagrange and
believed that he had brought mathematics to an almost completed
state\ldots{} In this conception questions about the status of
infinitesimals or the meaning of imaginary numbers were very much
secondary, and the primary emphasis was on operations, functions,
relations and the active process of working to solutions.
\cite[p.\;37]{Fr15}
\end{quote}

While initially plausible, this thesis dissolves under closer scrutiny
of the actual details.  The plausibility of Fraser's thesis hinges on
equivocation on the meaning of the term \emph{algebraic} used by
Fraser.  What exactly is meant by Lagrange's \emph{algebraic}
approach?  One can envision the following three possibilities for the
meaning of this term:
\begin{enumerate}
\item
Lagrange's approach to function theory based on power series
expansion;
\item
the principle of the so-called \emph{generality of algebra}
extensively relied upon by Lagrange following his predecessors;
\item
the algebraic foundations for analysis as initiated in Euler's
\emph{Introductio} \cite{Eu48}.
\end{enumerate}
Item (1) will not do since Cauchy extensively deals with nonanalytic
functions.  Thus, (1) is not plausible as a general description of
19th century work, since it leaves out Cauchy, who certainly was an
important figure for 19th century analysis.

Item (2), rather than being the context of Cauchy's work, was on the
contrary the main source of his dissatisfaction with Lagrange, as
detailed in the introduction to Cauchy's \emph{Cours d'Analyse}
\cite{Ca21}.

This leaves us with (3) Euler's \emph{Introductio} and subsequent work
by Lagrange and others.  Fraser's assumption of a clean separation
between algebra and infinitesimals may be compatible with modern
usage, but it is an anachronistic way of reading Euler's work where
algebraic manipulations with infinitesimals were ubiquitous; for
details see Bascelli et al.\;\cite{16a}.  In fact, Lagrange himself
eventually came to embrace infinitesimals in the \emph{second} edition
of his \emph{M\'ecanique Analytique}; see Katz--Katz\;\cite{11b}.

Thus, Euler's work can be viewed as an algebraic approach dealing with
algebraic manipulations of infinitesimal \emph{and} non-infinitesimal
entities alike.  Hence (3) is not a compelling argument that
infinitesimals only played a minor role.  The actual details are not
kind to Fraser's thesis.

\subsection
{Fraser, hyperreal analysis, and the details}
\label{s52}

Formulating assertions concerning a mathematical theory would require
that a scholar possess a basic level of competence with regard to the
theory, particularly if he seeks to base sweeping historical
conclusions on such claims.  In this section we will examine the level
of Fraser's competence in the matter of the hyperreals.

In his discussion of the hyperreals, Fraser finally provides some
details while attempting to summarize Robinson's construction of the
hyperreals via the compactness theorem applied to a suitable
language~$K$.

What is alarming is that Fraser \cite[p.\;25]{Fr15} specifically
describes his collection~$K$ as consisting of sentences of the
form~$c<v$ for each~$c\in\R$ (and no other sentences) in order to
force an infinite~$v$ in the model.  Earlier in that paragraph, Fraser
does comment on a ``formalized language that is rich enough to
formulate all sentences that are true in the real numbers~$\R$,'' but
this comment is not used in any way in his definition of~$K$ (and
furthermore contains a crucial inaccuracy: transfer applies to
\emph{first-order} sentences only, not to \emph{all sentences} as
Fraser claims).  Now any proper ordered field extension of~$\R$ will
satisfy the sentences of Fraser's collection~$K$, for example the
field of rational functions ordered by the behavior of the function at
infinity (i.e., a suitable lexicographic order).

Neither the compactness theorem nor hifalutin' mathematical logic are
required to build a proper ordered field extension of~$\R$.  Thus
Fraser entirely misses the point of the hyperreals as the only
extension possessing a full transfer principle (see
Section~\ref{f10}).  The transfer principle connects procedurally to
Leibniz's \emph{Law of continuity}.  Leibniz's theoretical strategy in
dealing with infinitesimals was analyzed in Katz--Sherry\;\cite{13f},
Sherry--Katz \cite{14d}, Bascelli et al.\;\cite{16a}, B\l aszczyk et
al.\;\cite{17c}, as well as the study by Bair et al.\;\cite{13a} that
Fraser is reacting to.  Leibniz's theoretical strategy was arguably
more robust than George Berkeley's flawed critique thereof.

We therefore cannot agree with Fraser's claim to the effect that
Leibniz allegedly ``never developed a coherent theoretical strategy to
deal with [infinitesimals]'' in \cite[p.\;31]{Fr15}.  Fraser's
attitude toward Leibniz's theoretical strategy in dealing with
infinitesimals is surely a function of his dismissive attitude toward
B-track historiography in general, which would similarly account for a
failure to appreciate the import of Robinson's extension.

\subsection{Transfer principle}
\label{f10}

The \emph{transfer principle} is a type of theorem that, depending on
the context, asserts that rules, laws or procedures valid for a
certain number system, still apply (i.e., are ``transfered'') to an
extended number system.  Thus, the familiar
extension~$\Q\hookrightarrow\R$ preserves the property of being an
ordered field.  To give a negative example, the
extension~$\R\hookrightarrow\R\cup\{\pm\infty\}$ of the real numbers
to the so-called \emph{extended reals} does not preserve the property
of being an ordered field.  The hyperreal
extension~$\R\hookrightarrow\astr$ preserves \emph{all} first-order
properties, e.g., the formula~$\sin^2 x + \cos^2 x =1$ (valid for all
hyperreal~$x$, including infinitesimal and infinite values
of~$x\in\astr$).  A construction of such an extension
$\R\hookrightarrow\astr$ appears in Section~\ref{s10}.  For a more
detailed discussion, see Keisler's textbook \emph{Elementary
Calculus}\;\cite{Ke86}.

\subsection{Does calculus become algorithmic over the hyperreals?}

Contrary to Fraser's claim cited in Section~\ref{s16}, algorithms
indeed play a central role in analysis exploiting the
extension~$\R\hookrightarrow\astr$.  A typical example is the
definition of the differential ratio~$\frac{dy}{dx}$ for a plane curve
as
\begin{equation}
\label{e41}
\text{st}\left(\frac{\Delta y}{\Delta x}\right)
\end{equation}
where~$\Delta x$ is an infinitesimal~$x$-increment and~$\Delta y$ is
the corresponding change in the variable~$y$.  Here ``st'' is the
standard part function which rounds off each finite hyperreal number
to the nearest real number.  Thus, Robinson's framework precisely
provides an algorithm for computing the differential ratio, which is
procedurally similar to Leibniz's transcendental law of homogeneity;
see Katz--Sherry \cite{13f} and Sherry--Katz \cite{14d} for an
analysis of the primary sources in Leibniz.

This modern B-track procedure is closer to the historical methods of
the pioneers of the calculus than to the (mathematically equivalent)
modern A-track procedures.  This is because the latter involve a
non-constructive notion of limit, where the value~$L$ of the limit
cannot be given algorithmically by a formula comparable to
\eqref{e41}, but rather has to be given in advance, so that one can
then elaborate an \emph{Epsilontik} proof that~$L$ is the correct
value.  In this sense, the B-track procedure as summarized
in~\eqref{e41} provides a better proxy for Leibniz's algorithm than
the A-track procedure, contrary to Fraser's claim.  For an analysis of
further constructive and algorithmic aspects of the hyperreal
framework see Sanders \cite{Sa16}, \cite{Sa17a}. \cite{Sa17b}.

\subsection{Euler's procedures and their proxies}
\label{s54}

Another example illustrating the advantage of B-track over A-track
when seeking proxies for the procedures of mathematical analysis as it
was practiced historically is Euler's definition of the exponential
function.  Fraser repeatedly mentions \emph{non-Archimedean analysis}
and \emph{non-Archimedean fields} in his article, but it is not
entirely clear whether he understands that the hyperreals are not
merely another non-Archimedean field.  In any such field there are
infinitesimals and infinitely large numbers incorporated in the
structure of an ordered field.  However, over the hyperreals one has
further structures, such as the theory of infinite, or more precisely
hyperfinite, sums like~$a_1+\ldots+a_N$ where~$N$ is an infinite
hyperinteger, as well as infinite products, which provide proxies for
Euler's procedures where A-track proxies are not available; see Bair
et al.\;\cite{17a}.  Thus we can follow Euler's procedure in writing
the exponential function~$e^{kz}$ as
\[
\left(1+\frac{kz}{i}\right)^i,
\]
where~$i$ is infinite number, in \cite{Eu48} vol.\;I, p.\;93, namely
an infinite product.  We can also follow Euler when he develops such a
product into an infinite sum using the binomial formula (for an
infinite exponent!), where A-track historians see only `dreadful
foundations' (see Section~\ref{s9}).

Euler's application of the binomial formula with an infinite exponent
is rewritten by Ferraro \cite[p.\;48]{Fe04} in
modern~$\sum_{r=0}^\infty$ notation relying on a Weierstrassian notion
of limit (``for every $\varepsilon>0$ there is an $N$ such that if
$n>N$ then the partial sum satisfies
$|L-\sum_{i=0}^n\ldots|<\varepsilon$, etc.'')  Yet Euler's procedure
here admits a closer proxy in terms of a \emph{hyperfinite} sum with
$i+1$ terms, where~$i$ is an infinite hyperinteger; see Bair et
al.\;\cite[p.\;222]{17a} for a discussion.  A related analysis of a
reductionist reading of Euler by H.\;M.\;Edwards appears in Kanovei et
al.\;\cite{15b}.

\subsection{Surreal blunder}
\label{s5}

The surreal number field was developed by Conway \cite{Co76} or
alternatively by Alling \cite{Al85}.  This field provides an example
of a B-continuum that possesses no transfer principle.%
\footnote{See Section~\ref{f10} for a brief introduction to the
transfer principle.}

Conway's maximal class surreal field is denoted \emph{No}.  Ehrlich
\cite{Eh12} refers to it as as the \emph{absolute arithmetic
continuum}.  This term has the correct connotations since the surreal
system is of some pertinence specifically in \emph{arithmetic}, being
the \emph{absolutely} maximal number system one could possibly
develop.  This may one day make the \emph{Guinness Book of Records},
but the field \emph{No} is only marginally relevant in analysis, since
even a simple function like the sine function cannot be extended to
all of \emph{No}.  Conway's system therefore cannot serve as a
foundation for analysis since \emph{No} does not satisfy transfer
beyond being a real closed field.

Namely, there is no surreally defined predicate of being an integer
which still honors transfer.  Conway's own version of surintegers fail
to satisfy transfer since~$\sqrt{2}$ turns out to be surrational.

On the other hand, Alling's surreals are field-isomorphic to a
properly defined hyperreal field with whatever transfer one wishes.
Thus the surreals do admit extensions of the usual relations of being
natural, etc., and the usual functions like the sine function, but
only via the hyperreals, rather than by means of the
\emph{sur}-construction.

Meanwhile, \cite[p.\;26]{Fr15} lists the surreals as one of the
possible frameworks for non-Archimedean analysis, alongside the
hyperreal framework and Bell's framework for Smooth Infinitesimal
Analysis.  But, as noted by Ehrlich \cite[Theorem~20]{Eh12}, the field
\emph{No} is isomorphic to a maximal (class) hyperreal field, and
therefore cannot be said to be an independent framework, especially
since any transfer principle in the surreals necessarily derives from
the hyperreals.  Including the surreals on par with these viable
frameworks for analysis is a misjudgment that puts into question
Fraser's technical ability to evaluate the appropriateness of modern
mathematical frameworks for interpreting the \emph{procedures} of the
mathematics of the past.

\section{Demonisation of Laugwitz and Euler}
\label{s6}

In this section we analyze Fraser's attitude toward the scholarship of
Detlef Laugwitz.  It turns out that Fraser repeatedly misrepresents
Laugwitz's work.  Fraser did this on at least two occasions: in his
2008 article for the \emph{New Dictionary of Scientific Biography},
and in his latest piece in 2015 for the monograph \emph{Delicate
Balance}.  The latter term can only be applied with difficulty to his
one-sided presentation of the state of Cauchy scholarship today.

\subsection
{Laugwitz's scholarship on analysis in Cauchy}
\label{s61}

In the abstract of his 1987 article in \emph{Historia Mathematica},
Laugwitz is careful to note that he interprets Cauchy's sum theorem
``with his [i.e., Cauchy's] own concepts'':
\begin{quote}
It is shown that the famous so-called errors of Cauchy are correct
theorems when interpreted with his own concepts. \cite[p.\;258]{La87}
\end{quote}
In the same abstract, Laugwitz goes on to stress that
\begin{quote}
\emph{No assumptions} on uniformity or on nonstandard numbers are
needed. (emphasis added)
\end{quote}
Laugwitz proceeds to give a lucid discussion of the sum theorem in
terms of Cauchy's infinitesimals in section 7 on pages 264--266, with
not a whiff of modern number systems.  In particular this section does
not mention the article Schmieden--Laugwitz \cite{SL58}.  In a final
section 15 entitled ``Attempts toward theories of infinitesimals,"
Laugwitz presents a general discussion of how one might formalize
Cauchyan infinitesimals in modern set-theoretic terms.  A reference to
the article by Schmieden and Laugwitz appears in this final section
only.  Thus, Laugwitz carefully distinguishes between his analysis of
Cauchy's \emph{procedures}, on the one hand, and the
\emph{ontological} issues of possible implementations of
infinitesimals in a set-theoretic context, on the other.

\subsection
{Fraser's assumptions}
\label{s62}

Alas, all of Laugwitz's precautions went for naught.  In 2008, he
became a target of damaging innuendo in the updated version of
\emph{The Dictionary of Scientific Biography}.  Here Fraser writes as
follows in his 2008 article on Cauchy:
\begin{quote}
Laugwitz's thesis is that certain of Cauchy's results that were
criticized by later mathematicians are in fact valid \emph{if one is
willing to accept certain assumptions} about Cauchy's understanding
and use of infinitesimals.  These assumptions reflect a theory of
analysis and infinitesimals that was worked out by Laugwitz and \ldots
Schmieden during the 1950s.  \cite[p.\;76]{Fr08} (emphasis added)
\end{quote}
Fraser claims that Laugwitz's interpretation of Cauchy depends on
\emph{assumptions} that reflect a modern theory of infinitesimals.
Fraser's indictment, based on the Omega-theory of Schmieden and
Laugwitz (see e.g., \cite{SL58}), is off the mark, as we showed in
Section~\ref{s61}.  In the intervening years Fraser has apparently not
bothered to read Laugwitz's article\;\cite{La87} either.  Indeed,
Fraser's verdict is unchanged seven years later in 2015:
\begin{quote}
Laugwitz, \ldots{} some two decades following the publication by
Schmieden and him of the~$\Omega$-calculus[,] commenced to publish a
series of articles arguing that their non-Archimedean formulation of
analysis is well suited to interpret Cauchy's results on series and
integrals. \cite[p.\;27]{Fr15}
\end{quote}
What Fraser fails to mention is that Laugwitz specifically and
explicitly separated his analysis of Cauchy's \emph{procedures} from
attempts to account \emph{ontologically} for Cauchy's infinitesimals
in modern terms, as we showed in Section~\ref{s61}.

Fraser's dual strategy in his article involves praising his allies and
undermining his opponents.  He seeks to present an allegedly united
front of received A-tracker scholarship on Cauchy, while distancing
himself from Boyer whose bowing down to the triumvirate%
\footnote{See note~\ref{f1} on page \pageref{f1} for an explanation of
the term.}
even Fraser finds it difficult to defend.  Fraser also adopts \emph{ad
hominem} arguments when it comes to dealing with his opponents such as
Laugwitz, by suggesting that their goal in interpreting history
involves an inappropriate interposition of a lens of their
professional preoccupations in non-Archimedean mathematics.

Fraser's first line of attack, involving an allegedly united
\emph{Epsilontik} front, is as questionable as the second.  Cracks in
triumvirate unity are ubiquitous.  Thus, Schubring speaks with
surprising frankness about his disagreement with Grabiner's ideas on
Cauchy--Weierstrass:
\begin{quote}
I am criticizing historiographical approaches like that of Judith
Grabiner where one sees epsilon-delta already realized in Cauchy
\cite[Section~3]{Sc15}
\end{quote}
(see also B\l aszczyk et al.\;\cite{17e}).  Fraser similarly overplays
his hand when he quotes Grattan-Guinness's distinction between history
and heritage in \cite{GG04}.  For it is precisely in his 2004 article
that Grattan-Guinness called explicitly for a re-evaluation of Cauchy,
by putting into question the \emph{Epsilontik} track, and warns
against reading Cauchy as if he had read Weierstrass already:
\begin{quote}
The (post-)Weierstrassian refinements have become standard fare, and
are incorporated into the \emph{heritage} of Cauchy; but it is mere
feedback-style ahistory to read Cauchy (and contemporaries such as
Bernard Bolzano) as if they had read Weierstrass already \ldots{} On
the contrary, their own pre-Weierstrassian muddles need historical
reconstruction, and clearly \cite[p.\;176]{GG04} (emphasis added).
\end{quote}
It is precisely Fraser's \emph{Cauchy--Weierstrass} line that is being
referred to as \emph{heritage} (rather than \emph{history}) here.
Grattan-Guinness also wrote that Cauchy's proof of the sum theorem is
difficult to interpret because it is stated in infinitesimal terms
(see Section~\ref{s35}), acknowledging the limitations of the A-track
approach when it comes to interpreting Cauchy's 1853 text \cite{Ca53}
on the summation of series.  Weierstrass' followers broke with
Cauchy's infinitesimal mathematics, but one will not discover this by
reading Fraser's comments on Cauchy.

A curious aspect of Fraser's text is his sweeping claim against the
relevance of modern non-Archimedean theories to the history of
analysis.  Fraser's misunderstandings arguably arise from failing to
distinguish between procedure and ontology.  He fails to cite some of
the key studies of Euler, such as Kanovei \cite{Ka88} and
McKinzie--Tuckey \cite{MT}; Laugwitz's article \cite{La87} is cited
but not discussed, and instead misrepresented by Fraser (see
Section~\ref{s62}); Laugwitz's article \cite{La89} appears in Fraser's
bibliography but is not mentioned in the body of Fraser's article.

\subsection{Gray on Euler's foundations}
\label{s9}

The relevance of modern B-track theoretical frameworks can be
established by the fact that the routine A-tracker claims that Euler's
foundations are allegedly \emph{dreadful} can today be discarded in
favor of a far greater appreciation of the coherence of his
infinitesimal techniques (see e.g., Section~\ref{s54}).  Thus, Euler's
foundations are described as \emph{dreadfully weak} by Jeremy Gray:
\begin{quote}
Euler's attempts at explaining the foundations of calculus in terms of
differentials, which are and are not zero, are dreadfully weak.
\cite[p.\;6]{Gr08}
\end{quote}
Such sweeping pronouncements come at a high price in anachronism.
Characteristically, Gray does not provide any justification for such a
claim.  Such were indeed the received views prior to the work of
Robinson, Laugwitz, Kanovei, and others.

\subsection{Shaky arguments}

In a similar vein, Judith Grabiner talks about \emph{shaky arguments}
and \emph{unerring intuition}:
\begin{quote}
\ldots{} eighteenth-century mathematicians had an almost
\emph{unerring intuition}.  Though they were not guided by rigorous
definitions, they nevertheless had a deep understanding of the
properties of the basic concepts of analysis. This conclusion is
supported by the fact that many apparently \emph{shaky
eighteenth-century arguments} can be salvaged, etc.
\cite[p.\;358]{Gr74} (emphasis added)
\end{quote}
This passage of Grabiner's sheds no light on how Euler and other 18th
century mathematicians may have hit upon such \emph{unerring
intuitions} to guide their work in infinitesimal analysis, in the
first place.  Short of an account of mathematical intuition,
\emph{unerring intuition} does almost no explanatory work.  It merely
tells us that Euler, Lagrange, and others had a special power.  To be
sure, there are theories of mathematical intuition; those of Kant and
Brouwer come to mind.  But neither of those theories is primarily
concerned with a power for making mathematical discoveries.

The coherence of the procedures of an 18th century master like Euler
is better understood in light of their modern proxies within
consistent theories of infinitesimals, making Euler's \emph{arguments}
appear less \emph{shaky} than those of Fraser, Grabiner, and Gray.

\subsection{Summary of conclusions}
\label{s8}

We have examined Fraser's stated position and analyzed the assumptions
underpinning his text \cite{Fr15}.  A summary of our conclusions
follows.

\subsubsection{}
Fraser claims that classical analysis, understood as modern
Weierstrassian analysis, is a ``primary point of reference" for the
historiography of analysis, for the stated reason that it ``developed
out of the older subject.''  However, Fraser's claim involves circular
reasoning, since his inference depends on the assumption of a
teleological view of the history of analysis (as having the
Weierstrassian framework as its inevitable goal), criticized by Ian
Hacking in terms of the dichotomy of the \emph{butterfly model} and
the \emph{Latin model} for the evolution of mathematics.

\subsubsection{}
By adhering to the teleological view, Fraser commits himself to
Boyer's position while paying lip service to a critique of Boyer's
subservience to ``the great triumvirate'' and arithmetisation.

\subsubsection{}
Fraser's claim that Leibniz allegedly ``never developed a coherent
theoretical strategy to deal with [infinitesimals]'' reveals his lack
of familiarity with recent literature published in \emph{Erkenntnis},
\emph{Studia Leibnitiana}, and \emph{HOPOS}, elaborating the details
of just such a strategy of Leibniz's; see Katz--Sherry \cite{13f},
Sherry--Katz \cite{14d}, Bascelli et al.\;\cite{16a}.

\subsubsection{}
Fraser claims that no ``significant contribution to history" can
result from applying modern theories to the history of analysis, but
he appears to make an exception for modern Archimedean theories as
developed by Weierstrass and his followers, as we showed in the
context of a thought experiment involving (non\mbox{-)} Archimedean
issues in Section~\ref{s13}.

\subsubsection{}
The Grabiner--Fraser so-called ``Cauchy--Weierstrass" foundation was
and remains a \emph{tale} (for details see Borovik--Katz \cite{12b},
Katz--Katz \cite{12d}, Bascelli et al.\;\cite{18b}), or more
specifically an A-tracker ideological commitment, rather than an
accurate historical category.

\subsubsection{}
Fraser's claim that ``Robinson and Lakatos were mistaken in their
assessment of Cauchy'' is a non-sequitur exacerbated by Fraser's
failure to make explicit the alleged mistakes of Robinson and Lakatos.

\subsubsection{}
Fraser wishes to place ``Cauchy within a given intellectual and
historical milieu'' but in point of fact yanks Cauchy right out of his
milieu and insert him in a 1870 milieu.

\subsubsection{}
Fraser cites Grattan-Guinness on \emph{mathematical history} versus
\emph{mathematical heritage}, but by Gratan-Guinness's own standard
with regard to Cauchy, it is Fraser himself who indulges in
mathematical heritage at the expense of history.

\subsubsection{}
Fraser's perhaps most fundamental philosophical blunder is his
failure to appreciate the importance of the distinction of
\emph{procedure} versus \emph{ontology} when it comes to interpreting
the mathematics of the past.

\section{Construction of the hyperreals}
\label{s10}

Since certain historians find hyperreal numbers to be abstruse (see
Section~\ref{s36}), we include a construction of a hyperreal field to
demonstrate that it requires no more background than a serious
undergraduate course in algebra including the theorem on the existence
of a maximal ideal.

In an approach to analysis exploiting Robinson's framework, one works
with the pair~$\R\subseteq\astr$ where~$\R$ is the usual ordered
complete \emph{Archimedean} continuum, whereas~$\astr$ is a proper
extension thereof.  Such a field~$\astr$ could be called a
\emph{Bernoullian} continuum, in honor of Johann Bernoulli who was the
first systematically to use an infinitesimal-enriched continuum as the
foundation for analysis.  Robinson's field~$\astr$ obeys the transfer
principle (see Section~\ref{s52}).

A field~$\astr$ can be constructed from~$\R$ using sequences of real
numbers.  To motivate the construction, it is helpfpul to analyze
first the construction of~$\R$ itself using sequences of rational
numbers.  Let~$\Q^\N_C$ denote the ring of Cauchy sequences of
rational numbers.  Then
\begin{equation}
\label{51}
\R=\Q^\N_C\,/\,\text{MAX}
\end{equation}
where ``MAX'' is the maximal ideal in~$\Q^\N_C$ consisting of all null
sequences (i.e., sequences tending to zero).

The construction of a Bernoullian field can be viewed as a relaxing,
or refining, of the construction of the reals via Cauchy sequences of
rationals.  This can be motivated by a discussion of rates of
convergence as follows.  In the above construction, a real number~$u$
is represented by a Cauchy sequence~$\langle u_n\colon n\in\N\rangle$
of rationals.  But the passage from~$\langle u_n\rangle$ to~$u$ in
this construction sacrifices too much information.  We seek to retain
a bit of the information about the sequence, such as its ``speed of
convergence."  This is what one means by ``relaxing" or ``refining"
the equivalence relation in the construction of the reals from
sequences of rationals.

When such an additional piece of information is retained, two
different sequences, say~$\langle u_n\rangle$ and~$\langle
u'_n\rangle$, may both converge to~$u\in\R$, but at different speeds.
The corresponding ``numbers" will differ from~$u$ by distinct
infinitesimals.  For example, if~$\langle u_n\rangle$ converges to~$u$
faster than~$\langle u'_n\rangle$, then the corresponding
infinitesimal will be smaller.  The retaining of such additional
information allows one to distinguish between the equivalence class
of~$\langle u_n\rangle$ and that of~$\langle u'_n\rangle$ and
therefore obtain distinct hyperreals infinitely close to~$u$.  Thus,
the sequence~$\langle\frac{1}{n^2}\rangle$ generates a smaller
infinitesimal than~$\langle\frac{1}{n}\rangle$.

A formal implementation of the ideas sketched above is as follows.
Let us outline a construction of a hyperreal field~$\astr$.
Let~$\R^{\N}$ denote the ring of sequences of real numbers, with
arithmetic operations defined termwise.  Then we have
\begin{equation}
\label{52}
\astr=\R^{\N}/\,\text{MAX}
\end{equation}
where ``MAX'' is a suitable maximal ideal of the ring~$\R^\N$.  What
we wish to emphasize is the formal analogy between \eqref{51} and
\eqref{52}.%
\footnote{In both cases, the subfield is embedded in the superfield by
means of constant sequences.}

We now describe a construction of such a maximal ideal
$\text{MAX}=\text{MAX}_\xi$ in terms of a finitely additive measure
$\xi$.  The ideal MAX consists of all ``negligible''
sequences~$\langle u_n\rangle$, i.e., sequences which vanish for a set
of indices of full measure, namely,
\[
\xi\big(\{n\in\N\colon u_n=0\}\big)=1.
\]
Here~$\xi\colon\mathcal{P}(\N)\to\{0,1\}$ (thus~$\xi$ takes only two
values,~$0$ and~$1$) is a finitely additive measure taking the
value~$1$ on each cofinite set,%
\footnote{For each pair of complementary \emph{infinite} subsets
of~$\N$, such a measure~$\xi$ ``decides'' in a coherent way which one
is ``negligible'' (i.e., of measure~$0$) and which is ``dominant''
(measure~$1$).}
where~$\mathcal{P}(\N)$ is the set of subsets of~$\N$.  The
subset~$\mathcal{F}_\xi\subseteq\mathcal{P}(\N)$ consisting of sets of
full measure is called a \emph{free ultrafilter}.  These originate
with Tarski \cite{Ta30}.  The construction of a Bernoullian continuum
outlined above was therefore not available prior to that date.

The construction outlined above is known as an ultrapower
construction.  The first such construction appeared in \cite{He48}, as
did the term \emph{hyper-real}.  The transfer principle%
\footnote{See Section~\ref{f10} on page~\pageref{f10} for an
explanation of the term.}
for this extension is an immediate consequence of the theorem of
\Los{} \cite{Lo55}.  For a recent application in differential
geometry, see Nowik--Katz \cite{15d}.  For a survey of approaches to
Robinson's framework, see Fletcher et al.\;\cite{17f}.

\section{Conclusion}

We have exploited the dichotomy of A-track versus B-track methods in
the historical development of infinitesimal calculus in analyzing the
work of the great historical authors, including Leibniz, Euler, and
Cauchy.  We have paid attention also to the dichotomy of
\emph{practice} versus \emph{ontology}.  The historiography of
mathematics as practiced by historians whose training was limited to a
Weierstrassian framework falls short of the target when analyzing many
of the great historical authors, due to its failure to attend to the
dichotomies mentioned, as exemplified by work of L\"utzen and Fraser.

While we are sympathetic to Barabashev's sentiment that ``as a result
of the elaboration of non-standard analysis, Leibniz's discovery is
differentiated more and more from Newton's theory of fluxions and
fluents'' \cite[p.\;38]{Ba97}, we feel that a more convincing case
involves the \emph{relative} advantages of Robinson's framework in
interpreting the historical infinitesimal analysis as compared to a
Weierstrassian framework.

We have followed Ian Hacking's proposal to approach the history of
mathematics in terms of a Latin (nondeterministic) model of
development rather than a butterfly (deterministic) model.  Carl
Boyer's account of the historical development of the calculus is a
classic example of a deterministic model.  While in preliminary
declarations, Fraser seeks to distance himself from Boyer's view of
the history of analysis as inevitable progress toward
infinitesimal-\emph{frei} arithmetisation, actually, as we show, he
commits himself to a position similar to Boyer's.

We have argued in favor of applying modern mathematical techniques
beyond Weierstrass, specifically, foundational studies both in
geometry and arithmetic of real numbers to analyze Euclid's
\emph{Elements} and 19th century constructions of real numbers.  We
have argued against Fraser's claim that analysis of historical texts
that relies on modern mathematics beyond Weierstrass cannot produce a
significant contribution to history.  Our methodological premise
involves applying Robinson's framework for analysis with
infinitesimals to analyzing 17-, 18- and 19th century mathematical
texts, following the pioneering work of Abraham Robinson and Detlef
Laugwitz.

\section*{Acknowledgments} 

We are grateful to V. de Risi for bringing to our attention the work
by Claude Richard \cite{Ri45}.  V. Kanovei was partially supported by
the RFBR grant no.\;17-01-00705.  M. Katz was partially supported by
the Israel Science Foundation grant no.\;1517/12.

\end{document}